\documentclass[10pt]{article}
 \usepackage[latin1]{inputenc}
 \usepackage[dvips]{graphicx}
 \usepackage{wrapfig}
 \usepackage{amsmath}
 \usepackage{amsthm}
 \usepackage{amsfonts}
 \usepackage{amssymb}
 \usepackage{layout}
 \usepackage{verbatim}
 \usepackage{alltt}
 \usepackage{XYpic}
 \input xy
 \xyoption{all}
 
\oddsidemargin - 1.5 pt
\newcounter{thm}[section]

\renewcommand{\thethm}{\arabic{section}.\arabic{thm}}

\newenvironment{thm}{\par\medskip\noindent\refstepcounter{thm}
\bgroup{\hspace*{-0.15 cm}\bf{Theorem}
\thethm.}\bgroup\it}{\egroup \egroup\par\medskip}

\newenvironment{lemma}{\par\medskip\noindent\refstepcounter{thm}
\bgroup{\hspace*{-0.15 cm}\bf{Lemma} \thethm.}\bgroup\it}{\egroup
\egroup\par\medskip}

\newenvironment{prop}{\par\medskip\noindent\refstepcounter{thm}
\bgroup{\hspace*{-0.15 cm}\bf{Proposition}
\thethm.}\bgroup\it}{\egroup \egroup\par\medskip}

\newenvironment{cor}{\par\medskip\noindent\refstepcounter{thm}
\bgroup{\hspace*{-0.15 cm}\bf{Corollary}
\thethm.}\bgroup\it}{\egroup \egroup\par\medskip}

\newenvironment{question}{\par\medskip\noindent\refstepcounter{thm}
\bgroup{\hspace*{-0.15 cm}\bf{Question}
\thethm.}\bgroup\it}{\egroup \egroup\par\medskip}

\def\F{{\mathbb F}}
\def\Z{{\mathbb Z}}
\def\Q{{\mathbb Q}}
\def\P{{\mathbb P}}
\def\c{{\mathbb C}}
\def\r{{\mathbb R}}

\def\O{{\mathcal O}}
\def\L{{\mathcal L}}
\def\R{{\mathcal R}}

\def\M{{\mathcal M}}

\def\p{{\mathcal P}}
\def\f{{\mathcal F}}
\def\N{{\mathcal N}}
\def\h{{\mathcal H}}

\def\X{{\overline X}}
\def\C{{\overline C}}

\def\Fs{{\overline \F}}
\def\b{{\overline b}}

\def\S{{\tilde S}}

\def\Fq{{\F_q}}
\def\Fsq{{\Fs_q}}

\begin{document}

  \author{Seyfi Türkelli}
  \title{Counting Multisections in Conic Bundles over\\ a Curve defined over $\Fq$}

  \maketitle
  
\footnotetext[1]{{\it 2000 Mathematics Subject Classification}. Primary ; Secondary }

\begin{abstract}
For a given conic bundle $X$ over a curve $C$ defined over $\Fq$, we count irreducible branch covers of $C$ in $X$ of degree $d$ and height $e\gg 1$. As a special case, we get the number of algebraic numbers of degree $d$ and height $e$ over the function field $\Fq (C)$.
\end{abstract}
  
\section{Introduction}

 Let $L$ be a number field of degree $n$ over $\Q$. One defines the counting function $\N_L(d,m,B)$ to be the cardinality of the set of $m$-tuples of algebraic numbers $\alpha=(\alpha_1,...,\alpha_m)$ such that $[L(\alpha_1,..,\alpha_m):L]=d$ and $\h(\alpha)\leq B$ where $\h$ is the absolute height function on $\overline{L}^m$. It is a nontrivial fact that this set is finite-- a result due to Northcott \cite{no}. Determining the asymptotics for the counting function $\N_L(d,m,B)$ is an old problem with few results in the literature.

 In 1979, Schanuel proved that $\N_L(1,m,B)\sim S_L(m)B^{n(m+1)}$ where $S_L(m)$ is a constant depending only on $L$ and $m$, see \cite{sc}. In 1995, Schmidt considered the quadratic case, that is $d=2$, and proved that $\N_{\Q}(2,m,B)$ is asymptotic to $S(d,m)B^{2(m+1)}$ if $m\geq 3$, to $S(d,m)B^6logB$ if $m=2$ and to $S(d,m)B^6$ if $m=1$. In 1995, Gao showed that $\N_{\Q}(d,m,B)\sim S(d,m)B^{d(m+1)}$ for $d\geq 3$ and $d\leq m$ where $S(d,m)$ is a positive constant depending only on $d$ and $m$, see \cite{ga}. For $d\geq 3$ and $d> m$, Gao also proved that the order of growth of $\N_{\Q}(d,m,B)$ is $B^{d(m+1)}$. 

 Another result along this direction is due to Masser and Vaaler, see \cite{ma}. They proved that $\N_L(d,1,B)\sim s(L,d)B^{nd(d+1)}$ where $n=[L:\Q]$ and $s(L,d)$ is a constant depending only on $L$ and $d$. Their idea for $L=\Q$ is as follows: every algebraic number $\alpha\in \overline{\Q}$ with $[\Q(\alpha):\Q]=d$ and $\h(\alpha)\leq B$ corresponds to a monic polynomial $f(x)\in \Z[x]$ of degree $d$ such that $\h(\alpha)=\M_{\infty}(f)^{1/d}$ where $\M_{\infty}$ is the Mahler measure. Now, one counts the monic polynomials $f(x)\in \Z[x]$ of degree $d$ with $\M_{\infty}(f) < B^d$. 

 One can formulate the problem above for $L=\Fq(C)$ where $C$ is a projective smooth curve defined over $\Fq$. In this setting, we use the additive height function $h$ on $\overline{L}$. The difficult part of Masser and Vaaler's proof is mainly the technicality due to Mahler measure. In the case $L=\Fq(C)$, the height function is a lot more canonical. In other words, there is no complications due to the place at infinity. For example, the relation between the additive height of an algebraic number $\alpha$ of degree $d$ over $L$ and the height of its monic polynomial $f(x)\in L[x]$ is simply given by the equation $dh(\alpha)=h(f)$ where $h(f)$ is the height of the coefficient vector of $f(x)$, see \cite[Lemma 6]{th2}. Therefore, it is often more natural to work over function fields to understand the nature of the problem, and then formulate it over number fields.

 In case $L=\Fq(C)$, Thunder proves the function field version of Masser-Vaaler's result. To be precise, let $\M_L(d,1,e)$ be the number of algebraic numbers of degree $d$ over $L$ and with absolute additive height exactly $e$. Then, one has $\M_L(d,1,e)\sim d a(d)q^{end(d+1)}$ where $[L:\Fq(t)]=n$ and $a(d)$ is a constant depending only on $L$ and $d$ \cite[Theorem 1]{th2}. Note that having additive height $e$ corresponds to having Weil height $q^e$ and so the order of growth in Thunder's formula exactly matches the one in Masser-Vaaler's formula. Counting algebraic numbers over $\Fq(C)$ of degree $d$ is as the same as counting irreducible branch covers of $C$ of degree $d$ inside the ruled surface $\pi:\P^1\times C\rightarrow C$, which one can realize as a trivial conic bundle. In this paper, we make a modest generalization of this result by counting covers in any given conic bundle over a curve $C/\Fq$. 

 Let $q$ be a power of an odd prime, $C$ be a geometrically connected smooth reduced projective curve defined over $\Fq$ of genus $g$ and let $\mathcal{L}$ be an invertible sheaf over $C$ of degree $l$. Let $\pi:X\rightarrow C$ be a geometrically irreducible projective conic bundle without non-reduced fibers such that $X/\Fq$ is smooth. Let $\iota:X\hookrightarrow \P^2\times_{\Fq} C$ be an embedding given by the equation $ax^2+by^2+cz^2=0$ where $a,b,c \in \Gamma(C,\mathcal{L})$ with $Supp(a)\cap Supp(b)\cap Supp(c)=\emptyset$. Finally, let $f:=pr_1\circ\iota: X\rightarrow \P^2$ be the composition of $\iota$ with the projection map $pr_1$. 

 A \emph{multisection} is a branch cover of $C$ defined over $\Fq$ contained in $X$. If $\pi|_{Y}:Y\rightarrow C$ is a multisection then its \emph{degree} is just the degree of the branch cover. The \emph{height} of the multisection $Y\rightarrow C$ is defined to be the number $e(Y)=Y.H$ where $H$ denotes a Weil divisor in the divisor class attached to the line bundle $f^*\O_{\P^2}(1)$. We say a multisection is of \emph{type} $(d,e)$ if it is of degree $d$ and height $e$.

 Likewise, we define the \emph{height of a divisor} $D$ on $X$ to be the intersection number $e(D):=D.H$. We call $D$ \emph{of type} $(d,e)$ if $D.F=d$ for any geometric fiber $F$ and $D.H=e$.  
 
 In this paper, we estimate the number of irreducible multisections of type $(d,e)$ for any given smooth conic bundle $\pi:X\rightarrow C$. More precisely, let $$\zeta_C(s):=\prod_{\text{ closed }P\in C}\frac{1}{1-q^{-s{\rm deg}P}}$$ be the Riemann-zeta function of $C/\Fq$. Let $\M(d,e)=\M_X(d,e)$ be the number of irreducible multisections of type $(d,e)$ in $X$. Then, there is a constant $K(d)$ depending only on the conic bundle $X$ and $d$ reflecting the contribution of the singular fibres to the leading coefficient of our estimate of $\M(d,e)$, which will be defined in section 3, such that we have:    

\begin{thm}\label{thm_main}
 For any even positive integer $d\in \Z$ there exists a constant $N$ depending only on the conic bundle and $d$ such that for every $e>N$, we have $$\M(d,e)=\frac{J(C) K(d)a(d)}{(q-1)\zeta_C(d+1)}\sqrt{q}^{(d+1)e}+ O(\sqrt{q}^{de})$$ where $a(d)=\sqrt{q}^{(d+1)(2-2g)-l([d/2+1]^2-1)}$ and $J(C)$ is the cardinality of Jacobian of $C/\Fq$.
\end{thm} 

 We get Thunder's result for even integers $d\in \Z$ \cite[Theorem 1]{th2} as an immediate corollary of Theorem \ref{thm_main}. Note that, after having the necessary notation in section $3$, we will show that the corollary below is really Thunder's result.  
 
\begin{cor}\label{cor_main}
 Let $X= \P^1\times C$ and $X\hookrightarrow \P^2\times C$ be an embedding via $\O((2,1))$. Then, we have: $$\M_X(d,e)=\frac{ J(C)a(d)}{(q-1)\zeta_C(d+1)}\sqrt{q}^{(d+1)e}+ O(\sqrt{q}^{de})$$ for sufficiently big $e\gg 1$ where $a(d)=q^{(d+1)(1-g)}$ and $J(C)$ is the cardinality of Jacobian of $C/\Fq$.
\end{cor}




 The outline of the paper is as follows. In section 2, we will compute $PicX$ and the dimension of the linear system $|D|$ for any given divisor $D\in PicX$ of type $(d,e)$. In section 3, we first will determine the proportion of the fiber-free effective divisors in a given linear system $|D|$ and we will prove our theorem by estimating to the number of irreducible ones. We conclude our paper with a discussion on our conjecture.


\section{Dimension of a Linear System}

 In this section, our main purpose is, for a given a divisor $D$ of type $(d,e)$ on $X$, to compute the dimension of the linear system $|D|$ and also to find the proportion of irreducible curves in the linear system $|D|$, see Corollary \ref{cor_dim} and Corollary \ref{cor_sur}. First, we will compute some intersection numbers and using Riemann-Roch, we will get the Euler characteristic $\chi(X,\O(D))$. Second, in Proposition \ref{prop_vanH1}, we will show the vanishing of the first cohomology group for a divisor with sufficiently big height and this will give the dimension of the linear system and the proportion we desire.  

 Once and for all, fix a finite field $\Fq$ with $char\Fq >2$. Fix a smooth projective reduced geometrically irreducible curve $C/\Fq$ and a conic bundle $\pi:X\rightarrow C$ defined (and smooth) over $\Fq$ with an embedding $\iota:X\hookrightarrow \P^2\times_{\Fq} C$ and without non-reduced fibers. Fix also an even positive integer $d\in \Z$ and set $d':=d/2$. 

 Let the conic bundle $\iota:X\hookrightarrow \P^2\times_{\Fq} C$ be given by the equation $ax^2+by^2+cz^2=0$ where $a,b,c \in \Gamma(C,\mathcal{L})$ with $Supp(a)\cap Supp(b)\cap Supp(c)=\emptyset$ and $\mathcal{L}$ is an invertible sheaf over $C$. From now on, $l$ will denote the degree of the sheaf $\mathcal{L}$ and $H$ will denote a Weil divisor attached to $f^* O_{\P^2}(1)$. For a variety $V$, in $Div(V)$, linear equivalence will be denoted by "$\sim$" and numerical equivalence will be denoted by "$\equiv$". Genus of the base curve will be denoted by $g=g(C)$ and the \emph{arithmetic genus} of $X$ by $p_a(X):=\chi(\O_X)-1$. 
 
 Given a geometric or a closed  point $P$ on $ C$, we will denote the fiber over $P$ by $X_P:=X\times_C P$. Note that $X_P$ is either a smooth conic or a union of two distinct projective lines in which cases we call, respectively, \emph{nonsingular} and \emph{singular}. 

 In this section, we will mostly work over the algebraic closure $\Fsq$ to ease the notation and we will carry the results of this section to over $\Fq$. For this section, set $E_1,...,E_k,E'_1,...,E'_k$ to be the components of singular geometric fibres. Set $\X=X\times_{\Fq}\Fsq$ and $\C=C\times_{\Fq} \Fsq$.

\begin{lemma}\label{lemma_intnum}
 For any geometric fibre $F$ of $\pi$ and for any $i\in\{1,..,k\}$, we have

 \begin{enumerate}
\item[i.] $F^2=E_i.F=E_i'.F=0$;

\item[ii.] $E_i.E_i'=1$ and $E_i^2=E_i'^2=-1$;

\item[iii.] $H.F=2$ and $H.E_i=H.E_i'=1$;

\item[iv.] $H^2=l$;

\item[v.] $F.K_X=-2$ and $E_i.K_X=E_i'.K_X=-1$ where $K_X$ is the canonical divisor of $X$.
 \end{enumerate} 

\end{lemma}  

\begin{proof}
 Clearly, we have $(i),(iii)$ and that $E_i.E'_i=1$. Since $(E_i+E'_i)^2=0$ and $E_i.E'_i=1$, we have $E_i^2=E_i'^2=-1$. Setting $y=z=0$ in the equation of the conic bundle, we get the self-intersection $H^2={\rm deg}div(a)={\rm deg}\mathcal{L}=l$ where $div(a)=\sum_{P\in C, a(P)=0} v_P(a)P$ is the divisor on $C$ corresponding to section $a$.

 Now, $(v)$ immediately follows from the adjunction formula and this completes the proof of the lemma.   
\end{proof}

 If the generic fiber $X_\eta$ is a trivial, then one can show that there exists a section $C_0$ and that $Pic\X$ is generated by $\pi^*Pic\C$, $-1$-curves $E_i$'s and the section $C_0$. Otherwise, one replaces the section $C_0$ with $H$. Using Lemma \ref{lemma_intnum}, one can see that $H- 2C_0\in \pi^*Pic\C$. Therefore, in any case, if $D$ is a divisor with even intersection number $D.F$ then $D$ can be written as a sum of fibers, $-1$-curves and $H$. This is what we show below (one can indeed show that $\X$ is a blow up of a ruled surface over finitely many points):     

\begin{prop}\label{prop_pic}
Let $D$ be a divisor in $Pic(\X)$ such that $D.F$ is even. Then, $D=aH + \sum_{i=1}^r c_iF_i + \sum_{i=1}^k b_iE_i$ where $F_i$'s are geometric fibers and $E_i$'s are exceptional divisors. Moreover, $\pi^*:Pic(C)\hookrightarrow Pic(X)$ is injective.
\end{prop}

\begin{proof}


 Let $D$ be such a divisor on $\X$ and $F$ be a geometric fibre. Since $D.F$ is even, there exists an $m$ such that $D=D_0+mH$ with $D_0.F=0$. Consider the divisor
$$D_n=D_0+ nF.$$
\noindent We have:
$$D_n^2=D_0^2 \text{\indent and \indent} D_n.K_X=D_0.K_X-2n.$$
\noindent Applying Riemann-Roch, we get
\begin{eqnarray*}
l(D_n)-s(D_n)+l(K_X-D_n)&=& \frac{1}{2}(D_0^2 - D_0.K_X + 2n)+1+p_a(X)\\
                        &=& n+\chi(D_0).
\end{eqnarray*}
\noindent As $n\rightarrow \infty$, $l(K_X-D_n)$ vanishes (because the intersection number $(K_X-D_n).H$ tends to $-\infty$) and this implies that the linear system $|D_n|$ is nonempty for some $n$.  

 Let $D'\in |D_n|$. $D'.F=D_n.F=0$ so every component of $D'$ is either a pull-back of a divisor of $\C$ or a component of a singular geometric fibre {\it i.e.} $D'=\sum_i F_i + \sum_j E_j + \sum_k E'_k$ in $Pic(\X)$. Thus $$D=D_0+mH=D_n-nF+mH=D'-nF+mH$$
\noindent in $Pic(\X)$ and we are done with the first claim. 

 Since $C$ and $X$ are projective, the natural maps $PicC\hookrightarrow Pic\C$ and $PicX\hookrightarrow Pic\X$ are injective, see \cite[Proposition A.2.2.10]{hi}. Injectivity of $\pi^*:Pic \C \rightarrow Pic \X$ is clear because $ Pic\C\hookrightarrow Pic(\P^2\times_{\Fsq} \C)$  and $Pic(\P^2\times_{\Fsq} \C)\hookrightarrow Pic\X$ are injective.
        
\end{proof}

 We will need the following lemmas to compute the Euler characteristic:

\begin{lemma}\label{lemma_eulerH}
 We have $\chi(H)=4g-4+2l$ where $\chi(H):=2g(H)-2$ denotes the Euler characteristic of $H$. Consequently, $H.K_X=4g-4+l$.  
\end{lemma}

\begin{proof}
 By Bertini's theorem (recall that we are working over $\Fsq$), one can can choose a divisor $Y$ in the divisor class of $H$ which  is a smooth curve and degree-$2$ branch cover of $\C$ in $\X$.  In fact, let $\alpha\in \Fsq$ and $\O_{\P^2}(1):z=\alpha x$. Replacing $z$ by $\alpha x$ in the equation of $X$, one gets $Y:(a+\alpha^2c)x^2 + by^2=0.$ Since $a,b,c\in \Gamma(C,\L)$ has no common support, one can choose $\alpha\in \Fsq$ such that $Supp(a+\alpha^2c)\cap Supp(b)=\emptyset$ and this implies that $Y$ is a double cover of $\C$.
 
 Clearly, the ramification divisor of the cover $\pi|_{Y}:Y\rightarrow \C$ has degree $2l$. Hurwitz theorem implies that $\chi(H)=\chi(Y)=4g-4+2l$. By the adjunction formula and Lemma \ref{lemma_intnum}, we get $H.K_X=4g-4+l$.
\noindent   
\end{proof}

\begin{lemma}\label{lemma_cancl}
We have $K_X\equiv -H + (2g-2+l)F$ where $F$ is a geometric fiber.
\end{lemma}

\begin{proof}
By Proposition \ref{prop_pic}, we have $K_X\sim aH+bF+\sum_{i=1}^k c_iE_i$. By (iii) and (v) of Lemma \ref{lemma_intnum}, $2a=F.K_X=-2$ and $a=-1$. This implies that $E_i.K_X=-1-c_i$; by Lemma \ref{lemma_intnum} (v), we get $c_i=0$ for all $i=1,...,k$. By Lemma \ref{lemma_intnum}, we have $H.K_X=2b-l$ and by Lemma \ref{lemma_eulerH}, $b=(2g-2+l)$. We are done. 
\end{proof}

\begin{prop}\label{prop_euler}
Let $D\equiv d'H+aF+\sum_{i=1}^k b_iE_i$ be a divisor on $X$ of type $(d,e)$. Then,
$$\chi(X,\O(D))=\frac{1}{2}[(d+1)(e+2-2g)-l((d/2+1)^2-1)-\sum_{i=1}^k b_i^2]$$
\noindent where $\chi(X,\O(D))$ is the Euler characteristic of the line bundle $\O(D)$. 
\end{prop}

\begin{proof}
Using Lemma \ref{lemma_intnum} and Lemma \ref{lemma_cancl}, we compute the intersection number $D.(D-K_X)$:
\begin{eqnarray*}
   &=& d'(d'+1)l+ 2d'(a+2-2g-l)+ d'\sum_{i=1}^k b_i + (d'+1)(2a+ \sum_{i=1}^k b_i) - \sum_{i=1}^k b_i^2
\\ &=& 2a(d+1) + (d+1)\sum_{i=1}^k b_i + (d^2/4-d/2)l + d(2-2g) - \sum_{i=1}^k b_i^2
\\ &=& (d+1)(e+2-2g)-l(d^2/4+d) - \sum_{i=1}^k b_i^2 + 2g-2.
\end{eqnarray*}
On the other hand, $X$ is birational to $\P^1\times_{\Fq} C$, see \cite[Theorem III.4]{be}. Since the arithmetic genus of a projective surface is a birational invariant \cite[Section V, Corollary 2.5]{ha}, $p_a(X)=-g$. Now, by Riemann-Roch and a simple computation, the claim of the proposition follows.
\end{proof}

  Recall that we want to count irreducible multisections in linear systems of divisors of type $(d,e)$. Not all divisors of type $(d,e)$ contain a multisection: 

\begin{prop}\label{prop_b}
Let $D\equiv d'H+ aF + \sum_{i=1}^k b_i E_i$ with $b_i\geq 0$ for all $i=1,..,k$. Suppose that there exists a multisection $Y\in |D|$. Then, $0\leq b_i\leq d'$ for all $i=1,...,k$. 
\end{prop}

\begin{proof}
Since $Y$ is a irreducible multisection, for each $i$ we have $Y.E_i=d'-b_i\geq 0$. Done.
\end{proof}

 By Proposition \ref{prop_pic}, a divisor $D$ on $X$ of type $(d,e)$ is numerically of the form $D\equiv d'H + aF + \sum_{i=1}^k b_i E_i$. One can also assume that $b_i\geq 0$ by replacing $E_i$'s  with $E'_i$'s  if necessary. From now on, by a \emph{divisor of type} $(d,e)$, we mean a divisor $D\equiv d'H+ aF + \sum_{i=1}^k b_i E_i$ of degree $d$ and height $e$ with the extra condition: $0\leq b_i\leq d'$ for all $i=1,...,k$.
 \ \\

 Above, we calculated $PicX$ and the Euler characteristic of a divisor $D$. In order to get the dimension of the linear system of a divisor, we need to show the vanishing of the first and the second cohomology groups (at least, for the divisors with sufficiently big height). First, we need a lemma: 

\begin{lemma}\label{lemma_isoH0}
$H^0(X,j_* \O_{X_P}\otimes_{\O_X} \O(D))\cong H^0(X_P,j^*\O(D))$.
\end{lemma}

\begin{proof}
We know that $$H^0(X,j_*\O_{X_P}\otimes_{\O_X}\O(D))\cong Hom_{\O_X}(\O_X,j_*\O_{X_P}\otimes_{\O_X}\O(D));$$ by the projection formula \cite{ha}[II, Exercise 5.1], $$\indent\indent\indent\indent\indent\cong Hom_{\O_X}(\O_X, j_*j^*\O(D));$$ by the adjunction property,  $$\indent\indent\indent\indent\indent\cong Hom_{\O_{X_P}}(j^*\O_X,j^*\O(D)).$$ Finally, by definition, $j^*\O_X\cong \O_{X_P}$  and we are done.
\end{proof}

\begin{prop}\label{prop_vanH1}
There exists an absolute positive constant $N\in \Z$ such that for any divisor $D\equiv d'H+ aF + \sum_{i=1}^k b_i E_i$ with $0\leq b_i\leq d'$ for $i=1,..,k$, $H^1(\X,\O(D))=H^2(\X,\O(D))=0$ whenever $e(D)\geq N$.
\end{prop}

\begin{proof}
 Let $\bar{b}=(b_1,...,b_k)$ be a $k$-tuple with $0\leq b_i\leq d'$ for $i=1,..,k$. Recall that $k$ is the number of singular fibres. Let $D_{\bar{b}}=d'H+ \sum_{i=1}^k b_iE_i$ and let $F_1,...,F_r$ be nonsingular geometric fibers. We first claim that there exists a constant $N_{\bar{b}}$ depending only on the conic bundle and the divisor $D_{\bar{b}}$ such that $H^1(\X,\O(D_{\bar{b}}+\sum_i a_iF_i))=0$ whenever $\sum_i a_i\geq N_{\bar{b}}$.

 Set $D=D_{\bar{b}}+\sum_{i=1}^k a_iF_i$. Let $P$ be a geometric point on $C$ and $j:X_P\hookrightarrow \X$ be the inclusion. Since $d=X_P.D>0$, $H^1(X_P,j^*\O(D))=0$. By exercise 11.8 of \cite[section III]{ha}, its higher direct image $\R^i\pi_*\O(D)=0$ for $i>0$ in a neigborhood of $P$ in $\C$. This holds for all geometric points so $\R^i\pi_*\O(D)=0$ for $i>0$ on $\C$. By exercise 8.1 of \cite[section III]{ha}, we have $$H^1(\X,\O(D))\cong H^1(\C,\pi_*\O(D)).$$ 

 On the other hand, by the projection formula, $$\pi_*\O(D)\cong \pi_*\O(D_{\bar{b}}) \otimes_{\O_C} \O(\sum_{i=1}^r a_iP_i)$$ where $F_i=X_{P_i}$ for some geometric points $P_1,...,P_r$ on $C$. Note that $\pi_*\O(D_{\bar{b}})$ is a coherent sheaf because $\O(D_{\bar{b}})$ is coherent and $\pi: X\rightarrow C$ is a projective morphism, see corollary 5.20 of \cite[III]{ha}.

 Now, by proposition 5.3 of \cite[III]{ha}, there exists a constant $N_{\bar{b}}$ depending only on the sheaf $\O(D_{\bar{b}})$ and on the curve $\C$ such that $$H^1(\C,\pi_*\O(D_{\bar{b}}) \otimes_{\O_C} \O(\sum_{i=1}^k a_iP_i))=0$$ whenever $\sum_{i=1}^r a_i \geq N_{\bar{b}}$. This completes the proof of the claim. 

 Let $N'=max\{N_{\bar{b}}| \bar{b}=(b_1,..,b_k) \text{ with }  0\leq b_i\leq d' \text{ for all }i \}$ and let $N= 2N' + d'e(H) + kd'$. We claim that $N$ is the desired constant. Say $D=D_{\bar{b}}+\sum_{i=1}^k a_iF_i$ where $\bar{b}$ is a tuple as above. Note that the constant $N_{\bar{b}}$ in the claim above does not depend on the choice of nonsingular fibres $F_i$'s. Note also that $a(D):=\sum_i a_i \geq N'$ whenever $e(D)\geq N$. Having noted these, it is now clear from the claim above that $H^1(\X,\O(D))=0$ if $e(D) \geq N$. 
 
 As for the vanishing of the second cohomology group $H^2(\X,\O(D))$, by Serre duality, $H^2(\X,\O(D))\cong H^0(\X,\O(K_X-D))$. By Lemma \ref{lemma_cancl}, $K_X-D\equiv -(d'+1)H+ cF + \sum_{i=1}^k c_iE_i$ for some constants $c,c_1,...,c_k$. Since $(K_X-D).F=-d-2$, $H^0(\X,\O(K_X-D))=0$ and we are done.
 
\end{proof}

As an immediate corollary, we get the main result of this section:

\begin{cor}\label{cor_dim}
$l(D)=\chi(X,\O(D))$ whenever $e(D)\geq N$.
\end{cor}

 Clearly, an effective divisor $Y$ on $X$ contains a fiber $X_P$, $P\in C$ a closed point, if a section corresponding to $Y$ vanishes on $X_P$. Next corollary gives us the proportion of effective divisors containing $X_P$ when ${\rm deg}P$ is "small" relative to the height of the divisor. Given a closed point $P\in C$, denote the irreducible components of the fiber $X_P$ by $E_P,E_P'$ (in case the fiber is already irreducible, then set the convention: $E_P=E_P'=X_P$). Let $$\f:=\{E_P,E_P'| P\in C \text{ is a closed point }\}$$ be the set of all the irreducible components of all the fibers. Since any finite subset $S$ of $\f$ can be viewed as a divisor on $X/\Fq$, it makes sense to talk about its height $e(S):=S.H$ (as defined above). Given a finite subset $S\subseteq \f$, we denote the union of the fibers in $S$ by $\S$, which is a closed subscheme of X.

\begin{cor}\label{cor_sur}
Let $D$ be a divisor on $X$. If $S$ is a finite subset of $\f$ with $e(S)\leq e(D)-N$, then the natural map $$\phi_S:H^0(X,\O(D))\rightarrow H^0(\S,j^*\O(D))$$ is surjective where $j:\S\hookrightarrow X$ is the inclusion.
\end{cor}

\begin{proof}
 We consider $S$ as a divisor on $X$ in the obvious way. The following short exact sequence $$0\rightarrow \O(D-S)\rightarrow \O(D)\rightarrow j_*\O_{\S}\otimes_{\O_{X}} \O(D)\rightarrow 0$$ with Lemma \ref{lemma_isoH0}, induces the long exact sequence $$0\rightarrow H^0(X,\O(D-S))\rightarrow H^0(X,\O(D))\rightarrow H^0(\S,j^*\O(D))\rightarrow H^1(X,\O(D-S)).$$
\noindent Now, the claim follows immediately. 
\end{proof}


\section{Proof of Theorem \ref{thm_main}}

 For this section, set $C_{<r}$ to be the finite closed subscheme of $C$ containing all the closed points $P\in C$ with ${\rm deg}P < r$. Likewise, set $C_{\geq r}=C\backslash C_{<r}$ to be the open subscheme. Let $C^1$ be the finite proper subscheme of $C$ containing all the closed points $P\in C$ whose fiber $X_P$ is irreducible but not geometrically irreducible, $C^2$ be the finite proper subscheme of $C$ consisting of all the closed points $P\in C$ whose fiber $X_P$ is reducible and let $C^0=C\backslash (C^1\cup C^2)$ be the open subscheme. Finally, let $N$ be a constant as in Proposition \ref{prop_vanH1}.

 Recall that $E_P, E'_P$ denote irreducible components of the fiber $X_P$ for every closed point $P\in C^2$ and that $E_P=E'_P=X_P$ if $P\in C\backslash C^2$. By Proposition \ref{prop_b}, with a right choice of connected component, every divisor $D$ on $X$ of type $(d,e)$ is numerically of the form: $D\equiv d'H + aX_P + \sum_{i=1}^k b_i E_{P_i}$ with $0\leq b_i\leq d'/{\rm deg}(P_i)$ for $i=1,...,n$ and for some (any) closed point $P\in C$. 

 By $e\gg 1$, we mean $e>M$ for some constant $M$ depending only on the conic bundle $\pi:X\rightarrow C$ and $d$. We write $f(x)=O(g(x))$ if  $f(x)<Mg(x)$ for some constant $M$ depending only on the conic bundle and $d$. Every constant stated in this section depends on the conic bundle or on the curve so this dependence will not be stated below.    
 
  In this section, we count irreducible multisections in a linear system of a given divisor $D$ on $X$ of type $(d,e)$. Here is the idea: effective divisors in the linear system $|D|$ which are not irreducible multisections are the ones whose support includes an irreducible component of a fiber or which can be written as a sum of two (or more) multisections of lower degree. At the end of the section, we will see that the asymptotic of reducible multi-sections of latter kind lies in the error term; this is easy! The difficult part is to determine the asymptotic of the "fiber-free" effective divisors.

 Say $D$ is a divisor on $C$ of type $(d,e)$ with sufficiently large height $e$. Given a closed point $P\in C$ of "small" degree, Corollary \ref{cor_sur} gives us the proportion of the effective divisors $Y\in |D|$ not including $E_P$, which is approximately $1-q^{-{\rm deg}P(d+1)}$ or $1-q^{-{\rm deg}P(d+2)}$, or $1-q^{-{\rm deg}P(d_P+1)}$ depending on the case where $d_P:=D.E_P/{\rm deg}P$. The main problem is that the conditions for a divisor to include $E_P$ in its support are not independent at different closed points $P\in C$. An idea used in the proof of Bertini's theorem over finite fields \cite{po}-- due to Poonen-- will help us to overcome this difficulty. In this section, our main purpose is to show that the terms in the asymptotic due to this dependence lies in the error term. 
 
 First, we need to introduce some notation. Let $D$ be a divisor on $X$ of type $(d,e)$. We define the counting function of the \emph{ fiber-free } effective divisors in the linear system $|D|$ as $$\M_f(D):=\sharp \{Y\in |D|: E_P\nsubseteq Y \text{ and }  E'_P\nsubseteq Y \text{ for each closed point } P\in C\}$$ and $\M_f(d,e)$ to be the sum of all $\M_f(D)$ over all divisors $D$ of type $(d,e)$ in $PicX$. Recall that $$\f=\{E_P,E'_P|P\in C \text{ is a closed point} \}$$ to be the set of irreducible components of the fibers. Given a finite subset $S$ of $\f$ with inclusion $j:\S\hookrightarrow X$, denote the proportion of the sections vanishing on $S$ by $$\p(D,S):=\sharp Ker\phi_S/q^{l(D)}$$ where $\phi_S:H^0(X,\O(D))\rightarrow H^0(\S,j^*\O(D))$ is the restriction. Note that $\p(D,\emptyset)=1$. 
 
\begin{lemma}\label{lemma_tri}
Let $D$ be a divisor on $X$ of type $(d,e)$. Then, $$\M_f(D)=\frac{q^{l(D)}}{q-1}\sum_{\substack{S\subset \f \\ e(S)\leq e}}(-1)^{\sharp S}(\p(D,S)-\frac{1}{q^{l(D)}}).$$ 
\end{lemma}

\begin{proof}
By definition, we have $$\M_f(D)=\frac{q^{l(D)}}{q-1}\sum_{\substack{S\subset \f \\ finite}}(-1)^{\sharp S}(\p(D,S)-\frac{1}{q^{l(D)}}).$$ If $e(S)>e$, then $\p(D,S)=1/q^{l(D)}$ since the height of a curve is not negative; hence, the desired equality.
\end{proof}

In this section, we will show that $\M_f(d,e)$ is asymptotic to the function introduced in section $1$ and that the difference $\M_f(d,e)-\M(d,e)$ lies in the error term. This amounts to approximating to the sum of the proportions $\p(D,S)$ over finite subsets $S$ of $\f$.
 
 One might start attempting to approximate $\p(D,S)$ by pretending that the conditions of vanishing on various $E_P$'s are all independent. This leads us to define proportions $$\p'(D,S):=\prod_{E_P\in S^0}q^{-{\rm deg}P(d+1)}\prod_{E_P\in S^1}q^{-{\rm deg}(d+2)}\prod_{E_P\in S^2}q^{-{\rm deg}(d_P+1)}$$ where $S^i=\{E_P\in S| P\in C^i\}$ for $i=0,1,2$ and $d_P=E_P.D/{\rm deg}P$. They are actually equal for subsets $S$ of small height:  
 
 \begin{lemma}\label{lemma_sml}
 Let $D$ be a divisor on $X$ of type $(d,e)$ and let $S$ be a finite subset of $\f$ with $e(S)\leq e-N$. Then, 
$\p(D,S)=\p'(D,S).$
\end{lemma}

\begin{proof}
 Clearly, $\sharp H^0(\S,j^*\O(D))=1/P'(D,S)$. Now, the claim follows from Corollary \ref{cor_sur}.
\end{proof}

 As for the subsets $S$ of $\f$ of "medium" height, the next proposition shows that the number $$\p(D,med):=\sum_{\substack{S\subset \f  \\ e-N\leq e(S)\leq e}}(-1)^{\sharp S}\p(D,S)$$ lies in the error term. First, we need a lemma:
 
 \begin{lemma}\label{lemma_mederr}
 Let $D$ be any divisor on $X$. Given a subset $S\subset \f$, set $S^0=\{E_P\in S| P\in C^0\}$. Then, for $m\gg 1$, we have $$\sum_{\substack{S\subset \f  \\ e(S)=m}}\p'(D,S^0)=O(q^{-m/2})$$ where the implicit constant in question is independent of $m$. 
 \end{lemma}
 
 \begin{proof}
   Denote the number of closed points in a finite scheme $V/\Fq$ by $|V|$. We claim that for any given positive integer $n\in \Z$, one can take $r\gg 1$ sufficiently big such that $|V| \leq r/n$ for any finite subscheme $V$ of $C$ of degree $r$. Let $n$ be given. Then, $\sharp C(\F_{q^{2n-1}})<cq^{2n-1}$ for some constant $c$ depending only on the curve $C$ \cite{lw}. Choose $r$ such that $r>2ncq^{2n-1}$. Then, $$|V_{<2n}| \leq |C_{<2n}| < cq^{2n-1} < \frac{r}{2n} .$$ On the other hand, $$2n.|V_{\geq 2n}| < \sum_{P\in V_{\geq 2n}}{\rm deg}P \leq r $$ and so $|V_{\geq 2n}|\leq \frac{r}{2n}$. This completes the proof of the claim. 
 
 Since there are only finitely many singular fibers, ${\rm deg}\pi(S^0) > r:=[m/2]-u$ for some fixed $u\in \Z$ (depending only on $X$) and for every such set $S$ with $e(S)=m$. By definition, $\p'(D,S^0)=\prod_{E_P\in S^0}q^{-{\rm deg}P(d+1)} < q^{-r(d+1)}$. Therefore, $$\sum_{\substack{S\subset \f  \\ e(S)=m}}\p'(D,S^0)\leq \sum_{\substack{S\subset \f  \\ e(S)=m}}q^{-r(d+1)} \leq \sharp\{S\subset \f | e(S)=m\}q^{-r(d+1)}.$$ 
 
 Again, because there only finitely many singular fibers, ${\rm deg} \pi (S)< r+w$ for some integer $w$ (depending only on the conic bundle $X$) and for any such set $S$. One can see that $$\sharp \{S\subset \f | e(S)=m \} < 3^s  \sharp \{V\subset C | {\rm deg}V=r+w\}$$ where $s$ is the number of singular fibers. One can also see that $$\sharp \{ V\subset C | {\rm deg}V=r+w \} < p(r+w)c^l q^{r+w}$$ where $p(r+w)$ is the number of partitions of $r+w$, $c$ is the constant defined in the claim above and $l = max \{|V| : V\subset C \text{ and }{\rm deg}V=r+w \}$. 

 Now, one can choose $m\gg 1$ such that $p(r+w)< q^{(r+w)/2}$ by \cite[Theorem 14.5]{ap} and $l< \frac{r+w}{2log_qc}$ by the claim above. Summing it up, we have $$\sum_{\substack{S\subset \f  \\ e(S)=m}}\p'(D,S^0) \leq 3^s q^{-(d+1)r}q^{(r+w)/2}q^{(r+w)/2}q^{r+w}$$ $$\indent\indent
 = O(q^{(1-d)r})\leq  O(q^{-r})=O(q^{-m/2})$$ since $d\geq 2$. This completes the proof.
 \end{proof}

 \begin{prop}\label{prop_mederr}
 For any divisor $D$ of type $(d,e)$ with $e\gg 1$, we have $$\p(D,med)=O(q^{-e/2}).$$
 \end{prop}
 
 \begin{proof}
   Let $D$ be such a divisor with $e(D)> N$. We first claim that there exists a constant $A$ depending on the conic bundle and $d$ such that $\p(D,S)\leq A.\p'(D,S)$ for every set $S\subset \f$ with $e-N\leq e(S)\leq e$. By Proposition \ref{prop_vanH1}, $H^1(X,\O(D))=0$ and so the exact sequence $$0\rightarrow \O(D-S)\rightarrow \O(D)\rightarrow j_*\O_{\S}\otimes_{\O_{X}} \O(D)\rightarrow 0$$ induces the long exact sequence (as in Corollary \ref{cor_sur}) $$0\rightarrow H^0(X,\O(D-S))\rightarrow H^0(X,\O(D))\rightarrow H^0(\S,j^*\O(D))\rightarrow H^1(X,\O(D-S))\rightarrow 0.$$ This implies that $$\p(D,S)=\frac{\sharp H^0(X,\O(D-S))}{\sharp H^0(X,\O(D))}=\frac{\sharp H^1(X,\O(D-S))}{\sharp H^0(\S,j^*\O(D))}=\sharp H^1(X,\O(D-S)). \p'(D,S).$$ On the other hand, if $e-N\leq e(S)\leq e$ then $0\leq e(D-S)\leq N$. Since there are finitely many divisors $D'$ on $X$ with $0\leq e(D')\leq N$ and $D'.F=d$, we can take $A=\max\{\sharp H^1(X,\O(D'))| 0\leq e(D')\leq N \text{ and } D'.F=d\}$ and complete the proof of the claim.
   
  Now, by the claim above, we have $$\sum_{\substack{S\subset \f  \\ e-N\leq e(S)\leq e}}(-1)^{\sharp S}\p(D,S)\leq A\sum_{n=0}^N\sum_{\substack{S\subset \f  \\ e(S)=e-n}}\p'(D,S) \leq A\sum_{n=0}^N\sum_{\substack{S\subset \f  \\ e(S)=e-n}}\p'(D,S^0)$$ where $S^0=\{E_P\in S | P \in C^0\}$. By Lemma \ref{lemma_mederr}, for $e\gg 1$ we have $$\sum_{\substack{S\subset \f  \\ e(S)=e-n}}\p'(D,S^0)=O(q^{-e/2})$$ for every $n=0,...,N$ and hence we are done. 

 \end{proof}
 
 Let $C^2=\{P_1,...,P_n\}$ and $\b=(b_1,...,b_n)=(b_P)_{P\in C^2}$ with $-d'/{\rm deg}P_i\leq b_i\leq d'/{\rm deg}P_i$ for each $i=1,...,n$. Let $E$ be a divisor on $C$ of degree $1$ and let $F:=X\times_C E$ be a divisor on $X$ of height $2$. Choose a set of connected components $E_P$ of the fibers $X_P$ for each closed point $P\in C^2$ and define $D_{\b}:= d'H + aF + \sum_{P\in C^2} b_P E_P$ with $D.H=e$. Then, each divisor $D$ of type $(d,e)$ is numerically of the form $D\equiv D_{\b}$ for some tuple $\b$ as above. Consider the following constant $$K(D):=K(\b)=\prod_{P\in C^1}\frac{1-q^{-{\rm deg}P(d+2)}}{1-q^{-{\rm deg}P(d+1)}} \prod_{P\in C^2} \frac{(1-q^{-{\rm deg}P(d_P+1)})(1-q^{-{\rm deg}P(d'_P+1)})}{1-q^{-{\rm deg}P(d+1)}}$$ where  $d_P=D.E_P/{\rm deg}P=d'-b_P$ and $d'_P=D.E'_P/{\rm deg}P=d'+ b_P$. This constant depends only on the numerical equivalence class of $D$; more precisely, it depends on the tuple $\b$. This is the constant which reflects the contribution of singular fibres to the leading coefficient of the counting function $\M_f(D)$ for each relevant divisor $D$. Their contribution to $\M_f(d,e)$ is measured by the following constant: $$K(d):= \sum_{\b} \frac{K(\b)}{\sqrt{q}^{\sum_{P\in C^2} (b_P^2{\rm deg}P)}}$$ where $$\sum_{\b}:=\sum_{b_1=[-d'/{\rm deg}P_1]}^{[d'/{\rm deg}P_1]}......\sum_{b_n=[-d'/{\rm deg}P_n]}^{[d'/{\rm deg}P_n]}.$$ Notice that this constant depends only on $d$. Here, $b_i$'s run from $[-d'/{\rm deg}P_i]$ to $[d'/{\rm deg}P_i]$ because, at the beginning, we make a choice of irreducible component $E_P$ for $D_{\b}$ and we need to take each choice into account. We need the following lemma:

\begin{lemma}\label{lemma_smlpr}
 For any divisor $D$ of type $(d,e)$ on $X$ with $e\gg 1$, we have $$\p(D,sml)=\frac{K(D)}{\zeta_C(d+1)} + O(q^{-e/2})$$ where $\p(D,sml):=\sum_{\substack{S\subset \f  \\ e(S)< e-N}}(-1)^{\sharp S}\p(D,S)$.
\end{lemma}

\begin{proof}
By Lemma \ref{lemma_sml}, $$\p(D,sml)=\sum_{\substack{S\subset \f  \\ finite}}(-1)^{\sharp S}\p'(D,S) - \sum_{\substack{S\subset \f  \\ finite \\e(S)\geq e-N}}(-1)^{\sharp S}\p'(D,S).$$ On the other hand, one can easily see that $$\sum_{\substack{S\subset \f  \\ finite}}(-1)^{\sharp S}\p'(D,S)=\frac{K(D)}{\zeta_C(d+1)}.$$ By Lemma \ref{lemma_mederr}, for $e\gg 1$ we have $$\sum_{\substack{S\subset \f  \\ e(S)\geq e-N}}(-1)^{\sharp S}\p'(D,S) \leq \sum_{\substack{S\subset \f  \\ e(S)\geq e-N}}\p'(D,S^0)= \sum_{m=e-N}^\infty O(q^{-m/2})=O(q^{-e/2})$$ and we are done. 
\end{proof}

 With a little bit extra work, we will prove the main theorem as a consequence of following proposition:
 
\begin{prop}\label{prop_main}
 For sufficiently big $e\gg 1$, we have $$\M_f(d,e)=\frac{ J(C)K(d)a(d)}{(q-1)\zeta_C(d+1)}\sqrt{q}^{(d+1)e}+ O(\sqrt{q}^{de})$$ where $a(d)=\sqrt{q}^{(d+1)(2-2g)-l([d/2+1]^2-1)}$ and $J(C)$ is the cardinality of the Jacobian of $C$.
\end{prop}

\begin{proof}
 Let $e\gg 1$ and $E$ be a divisor on $C$ of degree $1$. Let $F:=X\times_C E$, which is a divisor on $X$ of height $2$. Given a tuple $\b=(b_P)_{P\in C^2}$, let $D_{\b}:=d'H+\sum_{P\in C^2}b_P E_P+ a_{\b}F$ be the divisor in $PicX$ of type $(d,e)$. Notice that $a_{\b}$ is uniquely determined by $e$ and the tuple $\b$.

 By Proposition \ref{prop_mederr} and Lemma \ref{lemma_smlpr}, we have $$\M_f(D_{\b})=\frac{K(\b)}{(q-1)\zeta_C(d+1)}q^{l(D_{\b})} + O(q^{l(D_{\b})-e/2})- \frac{1}{q-1}\sum_{\substack{S\subset \f  \\ e(S)\leq e}}(-1)^{\sharp S}.$$ One can apply the argument in Lemma \ref{lemma_mederr} and show that for $e\gg 1$ $$\sum_{\substack{S\subset \f  \\ e(S)\leq e}}(-1)^{\sharp S}\leq \sum_{\substack{S\subset \f \\ e(S)= e}}1= O(q^e)\leq O(q^{l(D_{\b})-e/2}).$$  
Now, let's compute
\begin{eqnarray*}
\M_f(d,e) &=& \sum_{\bar{b}}\sum_{D'\in J(C)}\M_f(D_{\b} + \pi^*D')\\
          &=& \sum_{\b}J(C)[ \frac{K(\b)}{(q-1)\zeta_C(d+1)}q^{l(D_{\b})} + O(q^{l(D_{\b})-e/2})];\\       
\end{eqnarray*}
\noindent since this is a finite sum whose size depends only on the number closed points $P\in C^2$ and on $d$, we can take the error term out, and thus get: 
\begin{eqnarray*}
          &=& \frac{J(C)}{(q-1)\zeta_C(d+1)}\sum_{\b}[K(\b)q^{l(D_{\b})}] + O(q^{(d+1)e/2-e/2});
\end{eqnarray*}
\noindent Using Proposition \ref{prop_euler} and Corollary \ref{cor_dim}, one can easily see that
\begin{eqnarray*}
          &=& \frac{J(C)}{(q-1)\zeta_C(d+1)}\sum_{\b}[\frac{K(\b)}{\sqrt{q}^{\sum_{P\in C^2}(b_P^2{\rm deg}P)}}]a(d)\sqrt{q}^{(d+1)e} + O(\sqrt{q}^{de});\\
          &=& \frac{ J(C)K(d)a(d)}{(q-1)\zeta_C(d+1)}\sqrt{q}^{(d+1)e}+ O(\sqrt{q}^{de}).
\end{eqnarray*}
\end{proof}
 
 We are now ready to prove the main theorem.
\ \\
\begin{thm}\label{main_thm} Using the notation above, for sufficiently big $e$, we have $$\M(d,e)=\frac{J(C) K(d)a(d)}{(q-1)\zeta_C(d+1)}\sqrt{q}^{(d+1)e}+ O(\sqrt{q}^{de})$$ where $a(d)=\sqrt{q}^{(d+1)(2-2g)-l([d/2+1]^2-1)}$ and $J(C)$ is the cardinality of Jacobian of $C$.
\end{thm} 

\begin{proof} 
 By Proposition \ref{prop_main}, we just need to show: $$\M_f(d,e)-\M(d,e)=O(\sqrt{q}^{de}).$$ Note that for a divisor $D$ of type $(d,e)$, $q^{l(D)}=O(\sqrt{q}^{(d+1)e}).$ If $Y$ is an effective fiber-free divisor of type $(d,e)$ which is not a prime, then $Y=Y_1+Y_2$ for some effective divisors $Y_1$ of type of type $(d_1,e_1)$ and $Y_2$ of type $(d_2,e_2)$ with $e(Y)=e_1+e_2$; we can also assume that $d'_1:=d_1/2\leq [d/4].$ Therefore,

\begin{eqnarray*}
\M_f(d,e)-\M(d,e) &\leq & \sum_{i=1}^{[d/4]}\sum_{j=0}^e O(\sqrt{q}^{(2i+1)j+(d-2i+1)(e-j)})\\
                  &\leq & \sum_{i=1}^{[d/4]}\sum_{j=0}^e O(\sqrt{q}^{(d-2i+1)e})=\sum_{j=0}^e \sum_{i=1}^{[d/4]} O(\sqrt{q}^{(d-2i+1)e})\\
                  &\leq & \sum_{j=0}^e O(\sqrt{q}^{(d-1)e})=(e+1)O(\sqrt{q}^{(d-1)e})\leq O(\sqrt{q}^{de}).
\end{eqnarray*}
This completes the proof.
\end{proof}

 Let $X=\P^1\times_{\Fq} C$ with an embedding $X\hookrightarrow \P^2\times_{\Fq}C$ via $\O(2,1)$. Then, $l=0$ and $a(d)=q^{(d+1)(1-g)}$; $K(d)=1$ since there is no singular fibre. For a given multisection $Y$ in $X$ corresponding a degree-$d$ algebraic number $P\in \overline{\Fq(C)}$ over $\Fq(C)$, $e(Y)=2nd.h(P)$ where $h(P)$ is the height function used in Thunder's paper \cite{th2} and $n=[\Fq(C):\Fq(t)]$. Having noted our notation, one can check that we get Thunder's result for even integers $d\in \Z$ \cite[Theorem 1]{th2} as an immediate corollary of Theorem \ref{main_thm}:

\begin{cor} The notation is being as above, we have: $$\M_X(d,e)=\frac{ J(C)a(d)}{(q-1)\zeta_C(d+1)}\sqrt{q}^{(d+1)e}+ O(\sqrt{q}^{de}).$$
\end{cor}

 Note that there are multisections of odd degree only if the conic bundle has trivial generic fibre, that is when the generic fibre is isomorphic to $\P^1_{\Fq(C)}$. We remark that one can get the result for odd $d$ in exactly the same way as above. Everything works! This case is not done here because one needs to split each result of this paper into cases where $d$ is odd or even. Since the case of trivial conic bundle is already proven, we find unnecessary to make the paper longer than it is. 
 
 Now, we will consider the number field analogue of our result and conclude the discussion with a question. So, let $X\hookrightarrow \P^2$ be a conic defined by the equation $ax^2+by^2+cz^2=0$ for some pair-wise relatively prime $a,b,c\in \O_L$ where $L$ is a number field with $[L:\Q]=n$. Let $H_X$ be the height of the coefficient vector of the equation of $X$ in $\P^2$. Let $\N_X(d,B)$ be the number of $P\in X(\overline{L})$ with $[L(P):L]=d$ and $H(P)<B$. Set $r$ to be the number of real embeddings of $L$; $2s$ to be the number of complex embeddings; $\Delta_L$ to be the norm of the discriminant; $h_L$ to be the class number; $w_L$ to be the number of roots of unity; $R_L$ to be the regulator and finally set $\zeta_L$ to be the Riemann-Zeta function of $L$.

 By \cite[Theorem]{ma}, $\N_{\P^1}(d,B)$ has the leading coefficient $$S(L,d)=dV_{\r}(d)^r V_{\c}(d)^s (d+1)^{r+s-1} (\frac{2^r(2\pi)^s}{\sqrt{|\Delta_L|}})^{d+1} \frac{h_LR_L}{w_L\zeta_L(d+1)}$$ where the embedding $\P^1\hookrightarrow \P^2$ in question is via $\O_{\P^1}(1)$. We want to discuss how the leading coefficient in this case generalizes to more general case where we have any conic $X\subseteq \P^2/L$. We will do this using our generalization in function field case and the analogy between function fields and number fields.

 The contribution (to the leading coefficient) of places of function field $\F_q(C)$ is $1/\zeta_C(d+1)$ and it does not change when passing from $\P^1$ to a conic $X\hookrightarrow \P^2$; in number field case, this contribution is $V_{\r}(d)^r V_{\c}(d)^s (d+1)^{r+s-1}/\zeta_L(d+1)$ which one does not expect to change by analogy (for their definitions, see \cite{ma}). Constants $h_LR_L$ and $w_L$ are, respectively, direct analogies of $J(C)$ and $q-1$; which will stay the same. The term $(2^r(2\pi)^s / \sqrt{|\Delta_L|})^{d+1}$ is analogous to $a(d)$ so it should be changed to $(2^r(2\pi)^s/ \sqrt{|\Delta_L|})^{d+1}(H_X^{d^2/4+d})^{-1}$. Note that constant $d$ in the leading coefficient of $S(L,d)$ does not appear in the leading coefficient of $\M(d,e)$ because $\M(d,e)$ counts degree-$d$ polynomials not their roots.

 The only thing left is the contribution of the "singular fibres", which is denoted by $K(X,d)$ and it is defined as follows. We may assume $X$ is defined over the ring of integers $\O_L$. Set $L^1$ to be the finite set of primes $\wp$ of $L$ such that $X\otimes \O_L/\wp$ is irreducible but not geometrically irreducible, and $L^2$ to be the finite set of primes $\wp$ of $L$ such that $X\otimes \O_L/\wp$ is reducible. Let $L^2=\{\wp_1,...,\wp_n\}$. Given a tuple $\b=(b_1,...,b_n)=(b_{\wp})_{\wp\in L^2}$, define $$K(\b)=\prod_{\wp\in L^1}\frac{1-N\wp^{-(d+2)}}{1-N\wp^{-(d+1)}}\prod_{p\in L^2} \frac{(1-N\wp^{-(d'+b_{\wp}+1)})(1-N\wp^{-(d'-b_{\wp}+1)})}{1-N\wp^{-(d+1)}}.$$ Finally, define $$K(X,L,d):= \sum_{\b} \frac{K(\b)}{{\prod_{\wp\in L^2} \sqrt{N\wp}^{b_{\wp}^2}}}$$ where $$\sum_{\b}:=\sum_{b_1=[-d'/N\wp_1]}^{[d'/N\wp_1]}......\sum_{b_n=[-d'/N\wp_n]}^{[d'/N\wp_n]}.$$

\begin{question}
 Let the notation be as above. Is it true that $$\N_X(d,B)=S(X,L,d)B^{nd(d+1)/2}+\O(B^{nd^2/2})$$ where  $$S(X,L,d)=\frac{K(X,L,d)}{H_X^{d^2/4+d}}dV_{\r}(d)^r V_{\c}(d)^s(d+1)^{r+s-1}(\frac{2^r(2\pi)^s}{\sqrt{|\Delta_L|}})^{d+1}\frac{h_LR_L}{w_L\zeta_L(d+1)}.$$ 
\end{question}

\ \\

   \vspace{.2 in}
   \textsc{department of mathematics, university of wisconsin, 480 Lincoln Dr Madison wi 53706}
   
   \textit{E-mail address:} turkelli@math.wisc.edu


\begin{thebibliography}{100}
\bibitem[Ap]{ap} T. M. Apostol, Introduction to Analytic Number Theory, Springer-Verlag, New York, 1997.
\bibitem[Be]{be} A. Beauville, Complex Algebraic Surfaces, London Math. Society Student Texts 34, Cambridge Univ. Press, Cambridge, 1996.
\bibitem[Ga]{ga} X. Gao, On Northcott's Theorem, Ph.D Thesis, University of Colorado (1995).
\bibitem[Ha]{ha} R. Hartshorne. Algebraic Geometry. Number 52 in GTM. Springer-Verlag, 1977.
\bibitem[HS]{hi} M. Hindry and J. Silverman. Diophantine Geometry An Introduction, volume 201 of Graduate Texts in Mathematics. Springer-Verlag, 2000.
\bibitem[MV]{ma} D. Masser and J. Vaaler, Counting algebraic numbers with large height II, preprint (2004).
\bibitem[Nr]{no} D. G. Northcott, An inequality in the theory of arithmetic on algebraic varieties, Proc. Cambridge
Phil. Soc. 45 (1949), 502-509 and 510-518.
\bibitem[Po]{po} B. Poonen, Bertini theorems over finite fields, Ann. of Math. (2) 160 (2004), no. 3, 1099-1127.
\bibitem[Sc]{sc} S. Schanuel, Heights in number fields, Bull. Math. Soc. France 107 (1979), 433-449.
\bibitem[Th1]{th1} J.L. Thunder, Counting subspaces of given height defined over a function field, preprint (2006).
\bibitem[Th2]{th2} J.L. Thunder, More on Heights Defined over a Function Field, preprint 2006.
\bibitem[LW]{lw} S. Lang and A.Weil, Number of points of Varieties in Finite Fields, American J. Math. 76, (1954) 819-827.
\end{thebibliography}
\end{document}